# BOGOMOLOV CONJECTURE FOR CURVES OF GENUS 2 OVER FUNCTION FIELDS

ATSUSHI MORIWAKI

Abstract. In this note, we will show that Bogomolov conjecture holds for a non-isotrivial curve of genus 2 over a function field.

## 1. Introduction

Let $k$ be an algebraically closed field, $X$ a smooth projective surface over $k$, $Y$ a smooth projective curve over $k$, and $f : X \to Y$ a generically smooth semistable curve of genus $g \geq 2$ over $Y$. Let $K$ be the function field of $Y$, $\overline{K}$ the algebraic closure of $K$, and $C$ the generic fiber of $f$. Let $j : C(\overline{K}) \to \mathrm{Jac}(C)(\overline{K})$ be a morphism given by $j(x) = (2g-2)x - \omega_C$ and $\| \ \|_{NT}$ the semi-norm arising from the Neron-Tate height pairing on $\mathrm{Jac}(C)(\overline{K})$. We set

$$B_C(P; r) = \left\{ x \in C(\overline{K}) \mid \|j(x) - P\|_{NT} \leq r \right\}$$

for $P \in \mathrm{Jac}(C)(\overline{K})$ and $r \geq 0$, and

$$r_C(P) = \begin{cases} -\infty & \text{if } \#\left(B_C(P;0)\right) = \infty, \\ \\ \sup\left\{ r \geq 0 \mid \#\left(B_C(P;r)\right) < \infty \right\} & \text{otherwise.} \end{cases}$$

Bogomolov conjecture claims that, if $f$ is non-isotrivial, then $r_C(P)$ is positive for all $P \in \mathrm{Jac}(C)(\overline{K})$. Even to say that $r_C(P) \geq 0$ for all $P \in \mathrm{Jac}(C)(\overline{K})$ is non-trivial because it contains Manin-Mumford conjecture, which was proved by Raynaud. Further, it is well known that the above conjecture is equivalent to say the following.

**Conjecture 1.1** (Bogomolov conjecture)**.** If $f$ is non-isotrivial, then

$$\inf_{P \in \mathrm{Jac}(C)(\overline{K})} r_C(P) > 0.$$

Moreover, we can think the following effective version of Conjecture 1.1.

**Conjecture 1.2** (Effective Bogomolov conjecture)**.** In Conjecture 1.1, there is an effectively calculated positive number $r_0$ with

$$\inf_{P \in \mathrm{Jac}(C)(\overline{K})} r_C(P) \geq r_0.$$

*Date*: October, 1995 (Second version).
We would like to express our thanks to Prof. Liu for his nice remarks.

In [2], we proved that, if $f$ is non-isotrivial and the stable model of $f : X \to Y$ has only irreducible fibers, then Conjecture 1.2 holds. More precisely,

$$\inf_{P \in \mathrm{Jac}(C)(\overline{K})} r_C(P) \geq \begin{cases} \sqrt{12(g-1)} & \text{if } f \text{ is smooth,} \\ \sqrt{\dfrac{(g-1)^3}{3g(2g+1)}\delta} & \text{otherwise,} \end{cases}$$

where $\delta$ is the number of singularities in singular fibers of $f : X \to Y$. In this note, we would like to show the following result.

**Theorem 1.3.** *If $f$ is non-isotrivial and $g = 2$, then $f$ is not smooth and*

$$\inf_{P \in \mathrm{Jac}(C)(\overline{K})} r_C(P) \geq \sqrt{\frac{2}{135}\delta}.$$

## 2. Notations and ideas

In this section, we use the same notations as in §1. Let $\omega^a_{X/Y}$ be the dualizing sheaf in the sense of admissible pairing. (For details concerning admissible pairing, see [5] or [2].) First we note the following theorem. (cf. [5, Theorem 5.6] or [2, Corollary 2.3])

**Theorem 2.1.** *If $(\omega^a_{X/Y} \cdot \omega^a_{X/Y})_a > 0$, then*

$$\inf_{P \in \mathrm{Jac}(C)(\overline{K})} r_C(P) \geq \sqrt{(g-1)(\omega^a_{X/Y} \cdot \omega^a_{X/Y})_a},$$

*where $(\ \cdot\ )_a$ is the admissible pairing.*

From now on, we assume $g = 2$. By the above theorem, in order to get Theorem 1.3, we need to estimate $\left(\omega^a_{X/Y} \cdot \omega^a_{X/Y}\right)_a$. First of all, we can set

$$\left(\omega^a_{X/Y} \cdot \omega^a_{X/Y}\right)_a = \left(\omega_{X/Y} \cdot \omega_{X/Y}\right) - \sum_{y \in Y} e_y, \tag{2.2}$$

where $e_y$ is the number coming from the Green function of $f^{-1}(y)$. This number depends on the configuration of $f^{-1}(y)$. So, let us consider the classification of semistable curves of genus 2. Let $E$ be a semistable curve of genus 2 over $k$ and $E'$ the stable model of $E$, that is, $E'$ is a curve obtained by contracting all $(-2)$-curves in $E$. It is well known that there are 7-types of stable curves of genus 2. Thus, we have the classification of semistable curves of genus 2 according to type of $E'$ as in Table 1. (In Table 1, the symbol $A_n$ for a node means that the dual graph of $(-2)$-curves over the node is same as $A_n$ type graph.) The exact value of $e_y$ can be found in Table 2 and will be calculated in §3.

Next we need to think an estimation of $\left(\omega_{X/Y} \cdot \omega_{X/Y}\right)$ in terms of type of $f^{-1}(y)$. According to Ueno [4], there is the canonical section $s$ of

$$H^0(Y, \det(f_*(\omega_{X/Y}))^{10})$$

such that $d_y = \mathrm{ord}_y(s)$ for $y \in Y$ can be exactly calculated under the assumption that $\mathrm{char}(k) \neq 2, 3, 5$. The result can be found in Table 2. Prof. Liu points out that by works of T. Saito [3] and Q. Liu [1], the value $d_y$ in Table 2 still holds even if $\mathrm{char}(k) \leq 5$.

Let $\delta_y$ be the number of singularities in $f^{-1}(y)$. Then, by Noether formula,

$$\deg(\det(f_*(\omega_{X/Y}))) = \frac{\left(\omega_{X/Y} \cdot \omega_{X/Y}\right) + \sum_{y \in Y} \delta_y}{12}.$$

On the other hand, by the definition of $d_y$,

$$\sum_{y \in Y} d_y = 10 \deg(\det(f_*(\omega_{X/Y}))).$$

Thus, we have

$$\left(\omega_{X/Y} \cdot \omega_{X/Y}\right) = \sum_{y \in Y} \left(\frac{6}{5} d_y - \delta_y\right). \tag{2.3}$$

Hence, by (2.2) and (2.3),

$$\left(\omega^a_{X/Y} \cdot \omega^a_{X/Y}\right)_a = \sum_{y \in Y} \left(\frac{6}{5} d_y - \delta_y - e_y\right). \tag{2.4}$$

Therefore, using (2.4) and Table 2, we have the following theorem. (Note that an inequality

$$\frac{abc}{ab + bc + ca} \leq \frac{a + b + c}{9}$$

holds for all positive numbers $a, b, c$.)

**Theorem 2.5.** *If $f$ is non-isotrivial, then $f$ is not smooth and*

$$\left(\omega^a_{X/Y} \cdot \omega^a_{X/Y}\right)_a \geq \frac{2}{135} \delta,$$

*where $\delta = \sum_{y \in Y} \delta_y$.*

Note that non-smoothness of $f$ can be easily derived from the fact that the moduli space $M_2$ of curves of genus 2 is an affine variety.

## 3. Calculation of $e_y$

Let us start calculations of $e_y$. If the stable model of a fiber is irreducible, $e_y$ is calculated in [2]. Thus it is sufficient to calculate $e_y$ for II(a), IV(a,b), VI(a,b,c) and VII(a,b,c). In these cases, the stable model has two irreducible components. Let $f^{-1}(y) = C_1 + \cdots + C_n$ be the irreducible decomposition of $f^{-1}(y)$. We set

$$D_y = \sum_{i=1}^{n} (\omega_{X/Y} \cdot C_i) v_i,$$

where $v_i$ is the vertex in $G_y$ corresponding to $C_i$. Especially, we denote by $P$ and $Q$ corresponding vertexes to stable components. Then, $D_y = P + Q$. Let $\mu$ and $g_\mu$ be the measure and the Green function defined by $D_y$. In the same way as in the Proof of Theorem 5.1 in [2],

$$e_y = -g_\mu(D_y, D_y) + 4c(G_y, D_y),$$

where $c(G_y, D_y)$ is the constant coming from $g_\mu$. By the definition of $c(G_y, D_y)$,

$$c(G_y, D_y) = g_\mu(P, P) + g_\mu(P, D_y).$$

Therefore, we have

$$e_y = 7g_\mu(P,P) - g_\mu(Q,Q) + 2g_\mu(P,Q).$$

Here claim:

**Lemma 3.1.** *$g_\mu(P,P) = g_\mu(Q,Q)$. In particular,*

$$e_y = 6g_\mu(P,P) + 2g_\mu(P,Q).$$

*Proof.* By the definition of $c(G_y, D_y)$,

$$c(G_y, D_y) = g_\mu(P,P) + g_\mu(P, P+Q) = g_\mu(Q,Q) + g_\mu(Q, P+Q).$$

Thus, we can see $g_\mu(P,P) = g_\mu(Q,Q)$. □

In the following, we will calculate $e_y$ for each type II(a), IV(a,b), VI(a,b,c) and VII(a,b,c). First we present the dual graph of each type and then show its calculation.

**Type II(a).**

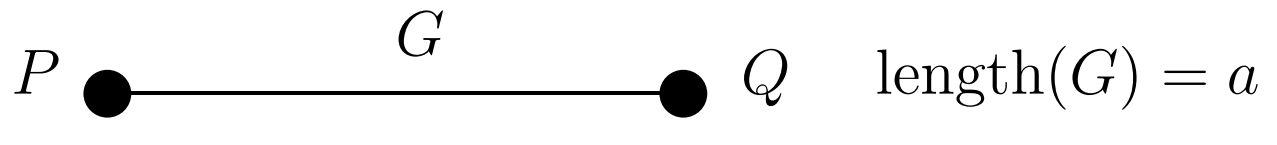

In this case, $\mu = \dfrac{\delta_P}{2} + \dfrac{\delta_Q}{2}$ by [5, Lemma 3.7]. We fix a coordinate $s : G \to [0, a]$ with $s(P) = 0$ and $s(Q) = a$. If we set

$$g(x) = -\frac{s(x)}{2} + \frac{a}{4},$$

then, $\Delta(g) = \delta_P - \mu$ and $\int_G g\mu = 0$. Thus, $g(x) = g_\mu(P, x)$. Hence

$$g_\mu(P,P) = \frac{a}{4} \qquad \text{and} \qquad g_\mu(P,Q) = -\frac{a}{4}.$$

Thus

$$e_y = 6g_\mu(P,P) + 2g_\mu(P,Q) = a.$$

**Type IV(a,b).**

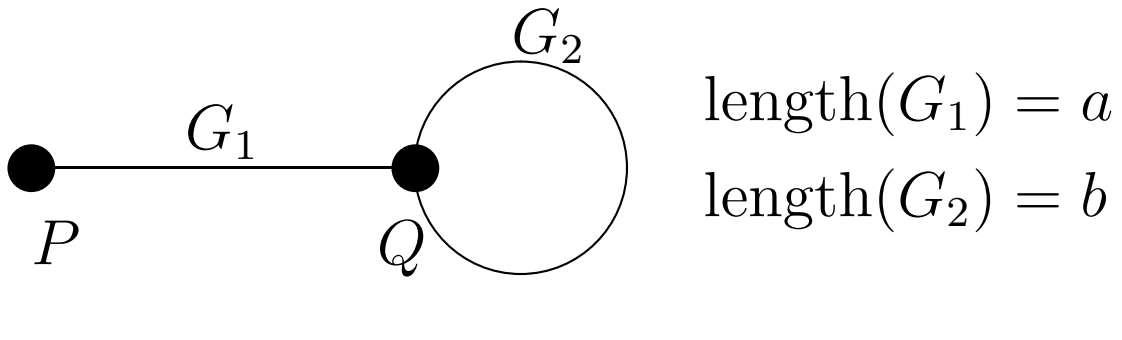

We fix coordinates $s : G_1 \to [0,a]$ and $t : G_2 \to [0,b)$ with $s(P) = 0$, $s(Q) = a$ and $t(Q) = 0$. In this case, $\mu = \dfrac{\delta_P}{2} + \dfrac{dt}{2b}$ by [5, Lemma 3.7]. We set

$$g(x) = \begin{cases} -\dfrac{s(x)}{2} + \dfrac{b+12a}{48} & \text{if } x \in G_1, \\[2ex] \dfrac{1}{2}\left(\dfrac{t(x)^2}{2b} - \dfrac{t(x)}{2}\right) + \dfrac{b-12a}{48} & \text{if } x \in G_2. \end{cases}$$

Then, $g$ is continuous, $\Delta(g|_{G_1}) = \dfrac{\delta_P}{2} - \dfrac{\delta_Q}{2}$, and $\Delta(g|_{G_2}) = \dfrac{\delta_Q}{2} - \dfrac{dt}{2b}$. Thus, $\Delta(g) = \delta_P - \mu$. Moreover, $\int_G g\mu = 0$. Therefore, $g(x) = g_\mu(P,x)$. Hence

$$g_\mu(P,P) = \frac{b+12a}{48} \qquad \text{and} \qquad g_\mu(P,Q) = \frac{b-12a}{48}.$$

Thus

$$e_y = 6g_\mu(P,P) + 2g_\mu(P,Q) = a + \frac{b}{6}.$$

**Type VI(a,b,c).**

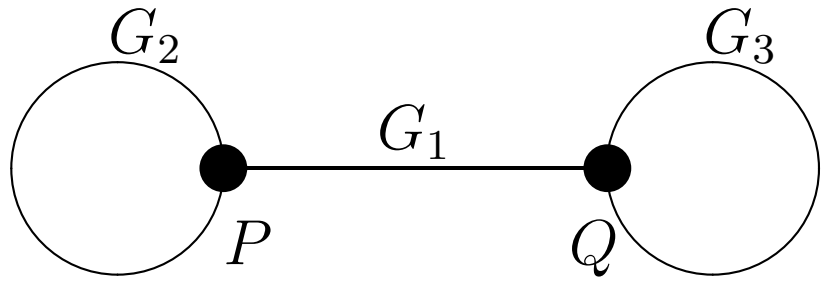

length$(G_1) = a$
length$(G_2) = b$
length$(G_3) = c$

We fix coordinates $s : G_1 \to [0,a]$, $t : G_2 \to [0,b)$ and $u : G_3 \to [0,c)$ with $s(P) = 0$, $s(Q) = a$, $t(P) = 0$ and $u(Q) = 0$. In this case, $\mu = \dfrac{dt}{2b} + \dfrac{du}{2c}$ by [5, Lemma 3.7]. We set

$$g(x) = \begin{cases} \dfrac{1}{2}\left(\dfrac{t(x)^2}{2b} - \dfrac{t(x)}{2}\right) + \dfrac{b+c+12a}{48} & \text{if } x \in G_2, \\[2ex] -\dfrac{s(x)}{2} + \dfrac{b+c+12a}{48} & \text{if } x \in G_1, \\[2ex] \dfrac{1}{2}\left(\dfrac{u(x)^2}{2c} - \dfrac{u(x)}{2}\right) + \dfrac{b+c-12a}{48} & \text{if } x \in G_3. \end{cases}$$

Then, $g$ is continuous, $\Delta(g|_{G_1}) = \dfrac{\delta_P}{2} - \dfrac{\delta_Q}{2}$, $\Delta(g|_{G_2}) = \dfrac{\delta_P}{2} - \dfrac{dt}{2b}$, and $\Delta(g|_{G_3}) = \dfrac{\delta_Q}{2} - \dfrac{du}{2c}$. Thus, $\Delta(g) = \delta_P - \mu$. Moreover, $\int_G g\mu = 0$. Therefore, $g(x) = g_\mu(P,x)$. Hence

$$g_\mu(P,P) = \frac{b+c+12a}{48} \qquad \text{and} \qquad g_\mu(P,Q) = \frac{b+c-12a}{48}.$$

Thus

$$e_y = 6g_\mu(P,P) + 2g_\mu(P,Q) = a + \frac{b+c}{6}.$$

**Type VII(a,b,c).**

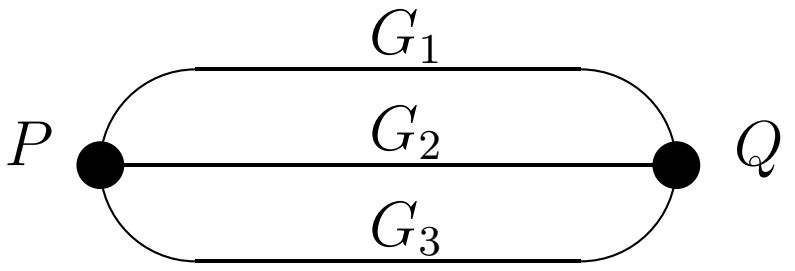

$\mathrm{length}(G_1) = a$
$\mathrm{length}(G_2) = b$
$\mathrm{length}(G_3) = c$

We fix coordinates $s : G_1 \to [0, a]$, $t : G_2 \to [0, b]$ and $u : G_3 \to [0, c]$ with $s(P) = 0$, $s(Q) = a$, $t(P) = 0$, $t(Q) = b$, $u(P) = 0$ and $u(Q) = c$. In this case, $\mu = \frac{ds}{3a} + \frac{dt}{3b} + \frac{du}{3c}$ by [5, Lemma 3.7]. We set

$$g(x) = \begin{cases} \dfrac{s(x)^2}{6a} - \left(\dfrac{1}{6} + \dfrac{1}{2}\dfrac{bc}{ab+bc+ca}\right)s(x) + \dfrac{a+b+c}{108} + \dfrac{1}{4}\dfrac{abc}{ab+bc+ca} & \text{if } x \in G_1, \\ \dfrac{t(x)^2}{6b} - \left(\dfrac{1}{6} + \dfrac{1}{2}\dfrac{ac}{ab+bc+ca}\right)t(x) + \dfrac{a+b+c}{108} + \dfrac{1}{4}\dfrac{abc}{ab+bc+ca} & \text{if } x \in G_2, \\ \dfrac{u(x)^2}{6c} - \left(\dfrac{1}{6} + \dfrac{1}{2}\dfrac{ab}{ab+bc+ca}\right)u(x) + \dfrac{a+b+c}{108} + \dfrac{1}{4}\dfrac{abc}{ab+bc+ca} & \text{if } x \in G_3. \end{cases}$$

Then, $g$ is continuous and $\Delta(g) = \Delta(g|_{G_1}) + \Delta(g|_{G_2}) + \Delta(g|_{G_3}) = \delta_P - \mu$. Moreover, $\int_G g\mu = 0$. Therefore, $g(x) = g_\mu(P, x)$. Hence

$$e_y = 6g_\mu(P, P) + 2g_\mu(P, Q) = \frac{2}{27}(a+b+c) + \frac{abc}{ab+bc+ca}.$$

Department of Mathematics, Faculty of Science, Kyoto University, Kyoto, 606-01, Japan
*E-mail address*: moriwaki@kusm.kyoto-u.ac.jp

TABLE 1. Classification of semistable curve $E$ of genus 2

| Type of $E$ | Description of the stable model $E'$ of $E$ | Figure of $E'$ and types of singularities by contracting $(-2)$-curves in $E$ |
|---|---|---|
| I | a smooth curve of genus 2 | 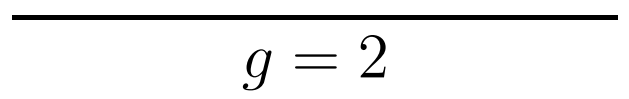  |
| II(a) | two elliptic curves meeting at one point transversally | 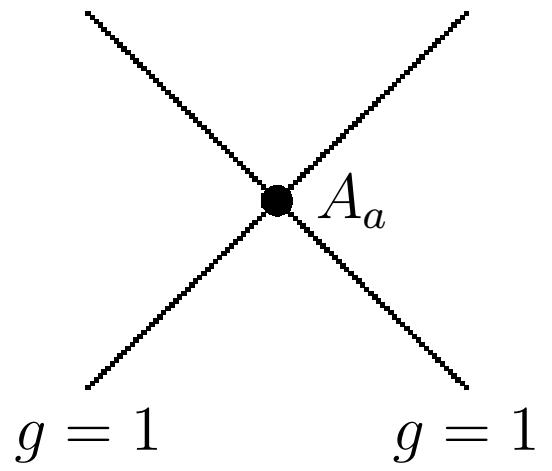  |
| III(a) | an elliptic curve with one node | 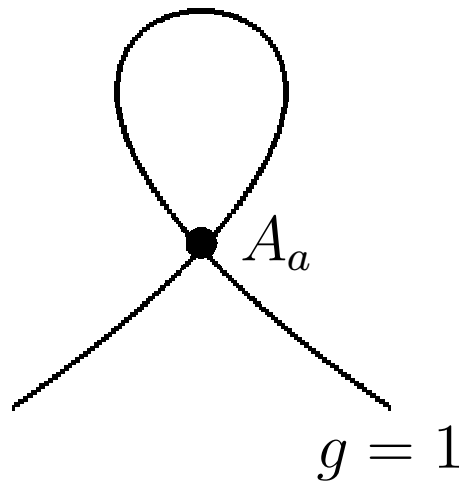  |
| IV(a,b) | a smooth elliptic curve and a rational curve with one node, which meet at one point transversally | 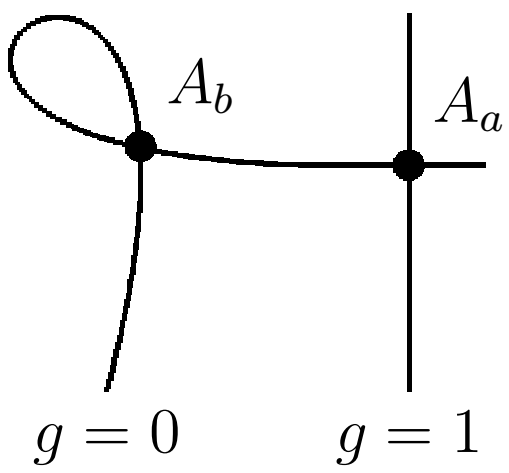  |
| V(a,b) | a rational curve with two nodes | 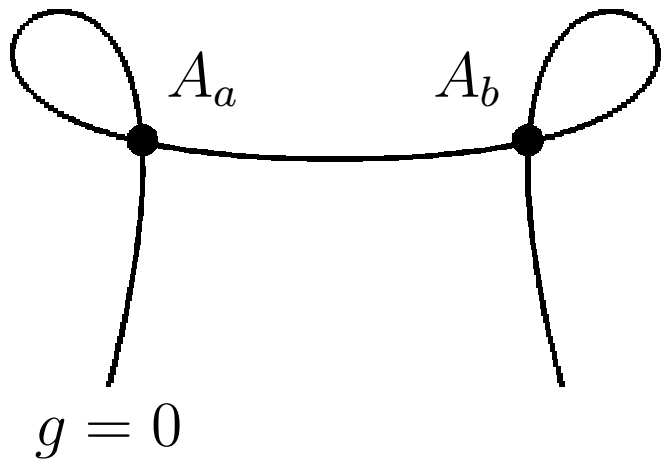  |

<table>
<tr><td>VI(a,b,c)</td><td>two rational curves with one node, which meet at one point transversally</td><td>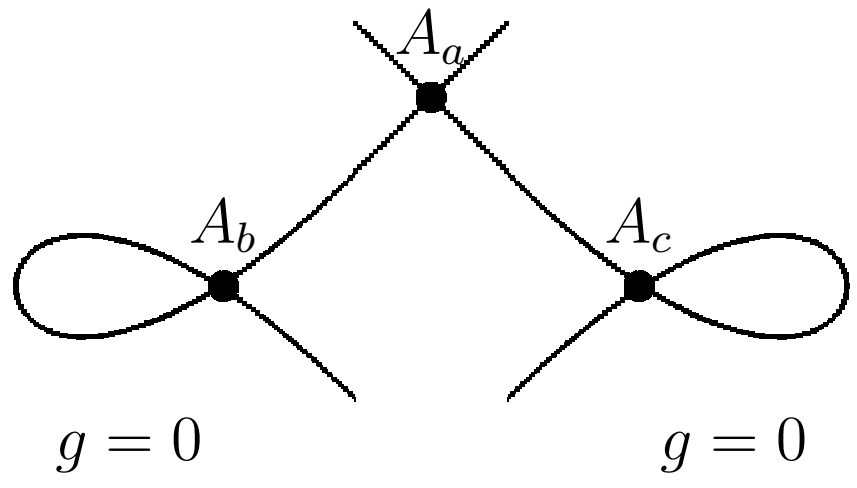
</td></tr>
<tr><td>VII(a,b,c)</td><td>two smooth rational curves, which meet at three points transversally</td><td>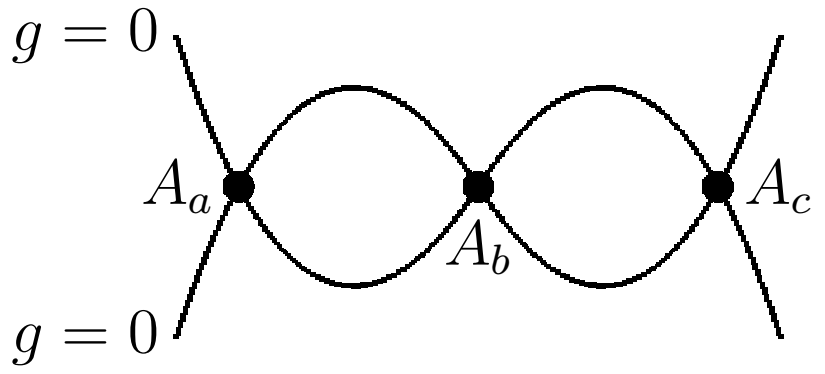
</td></tr>
</table>

TABLE 2. $\delta_y$, $d_y$ and $e_y$

| Type | $\delta_y$ | $d_y$ | $e_y$ |
|---|---|---|---|
| I | $0$ | $0$ | $0$ |
| II(a) | $a$ | $2a$ | $a$ |
| III(a) | $a$ | $a$ | $\frac{a}{6}$ |
| IV(a,b) | $a+b$ | $2a+b$ | $a+\frac{b}{6}$ |
| V(a, b) | $a+b$ | $a+b$ | $\frac{a+b}{6}$ |
| VI(a,b,c) | $a+b+c$ | $2a+b+c$ | $a+\frac{b+c}{6}$ |
| VII(a,b,c) | $a+b+c$ | $a+b+c$ | $\frac{2}{27}(a+b+c)+\frac{abc}{ab+bc+ca}$ |